\documentclass[a4paper,draft]{amsart}
\usepackage[a4paper,margin=25mm]{geometry}
\usepackage{amsmath,tikz}
\usepackage{amssymb}

\parskip\smallskipamount

\let\set\mathbb
\def\<#1>{\langle#1\rangle}

\def\Model#1#2#3#4#5{%
  \ref{#1} & $#4$\rule[-1em]{0pt}{2.5em} \\\hline
}
\def\nextcella{\global\let\nextcell\nextcellb&}
\def\nextcellb{\global\let\nextcell\nextcellc&}
\def\nextcellc{\global\let\nextcell\nextcelld&}
\def\nextcelld{\global\let\nextcell\nextcella\\\hline}
\global\let\nextcell\nextcella
\def\model#1#2#3#4#5{%
  \leavevmode\noindent
  \refstepcounter{model}\label{#1}%
  \rlap{\arabic{model}\mathstrut}%
  \kern1.15em
  \rlap{\raisebox{-15mm}{$#5,\dots\mathstrut$}}%
  \raisebox{-12mm}{\stepset{#4}{#2}}%
  \kern.6em
  \parbox[t]{.1\hsize}{\raggedright\footnotesize\leavevmode\let\comma\empty\writestepseta#2,@,.}%
  \nextcell
}
\newcounter{model}

\def\stepset#1#2{$\hbox{\begin{tikzpicture}[scale=.14]
    \ifx0#1\else\fill[lightgray](-4.5,-4.5) rectangle (4.5,4.5);\fi
    \foreach\x in {-4,...,4} \foreach\y in {-4,...,4} \draw(\x,\y) node {$\cdot$};
    \ifx0#1\draw[white] (0,0) node {$\bullet$};\else\draw[lightgray] (0,0) node {$\bullet$};\fi    
    \draw (-4.5,-4.5) rectangle (4.5,4.5)
    (-1.5,-4.5)--(-1.5,4.5) (1.5,-4.5)--(1.5,4.5) (-4.5,-1.5)--(4.5,-1.5) (-4.5,1.5)--(4.5,1.5);
    \drawstepset#2,@,\end{tikzpicture}}$}
\def\drawstepset#1,{\ifx@#1\let\next\empty\else\handlestep#1\let\next\drawstepset\fi\next}
\def\handlestep#1#2#3#4{%
  \def\x{0}\ifx#1-\def\x{-3}\fi\ifx#1+\def\x{3}\fi
  \def\y{0}\ifx#2-\def\y{-3}\fi\ifx#2+\def\y{3}\fi
  \ifx#3-\edef\x{\x-1}\fi\ifx#3+\edef\x{\x+1}\fi
  \ifx#4-\edef\y{\y-1}\fi\ifx#4+\edef\y{\y+1}\fi
  \draw(\x,\y) node {\tiny$\bullet$};
}
\def\writestepset#1,{\ifx@#1\let\next\empty\else\xhandlestep#1\let\next\writestepset\fi\next}
\def\u#1{\ifx#1-{\bar1}\fi\ifx#1+{1}\fi\ifx0#1{0}\fi}
\def\xhandlestep#1#2#3#4{\hbox{\footnotesize$\mathstrut\u{#1}\u{#2}\u{#3}\u{#4}$}}
\def\writestepseta#1,{\ifx@#1\let\next\empty\else\comma\def\comma{, }\xhandlestepa#1\let\next\writestepseta\fi\next}
\def\xhandlestepa#1#2#3#4{\hbox{$\mathstrut\u{#1}\u{#2}\u{#3}\u{#4}$}}

\newtheorem{tab}{Table}

\begin{document}

 \author[Manfred Buchacher]{Manfred Buchacher\,$^\ast$}
 \author[Sophie Hofmanninger]{Sophie Hofmanninger}
 \author[Manuel Kauers]{Manuel Kauers}
 \address{Manfred Buchacher, Institute for Algebra, J. Kepler University Linz, Austria}
 \email{manfred.buchacher@jku.at}
 \address{Sophie Hofmanninger, Institute for Algebra, J. Kepler University Linz, Austria}
 \email{sophie@hofmanninger.co.at}
 \address{Manuel Kauers, Institute for Algebra, J. Kepler University Linz, Austria}
 \email{manuel.kauers@jku.at}
 \thanks{$^\ast$ Supported by the Austrian Science Grant F5004.}

 \title{Walks with Small Steps in the 4D-Orthant}

 \begin{abstract}
   We provide some first experimental data about
   generating functions of restricted lattice walks
   with small steps in~$\set N^4$.
 \end{abstract}

 \maketitle

 \section{Introduction}

 For given sets $S\subseteq\{-1,0,1\}^2\setminus\{(0,0)\}$, many people have studied the number
 of lattice walks in $\set N^2$ starting at $(0,0)$ and consisting of $n$ steps each of
 which is taken from~$S$. If $a_{i,j,n}$ denotes the number of such walks ending at the
 point $(i,j)\in\set N^2$, then it is interesting to study the nature of the formal
 power series $f(x,y,t)=\sum_{n,i,j=0}^\infty a_{i,j,n}x^i y^j t^n$, i.e., whether it is
 algebraic, or, if not, whether it is at least D-finite. It turns out that the nature
 of $f$ depends on the choice of~$S$, and a lot of work has been done in order to
 classify all step sets $S$ according to the nature of the corresponding generating
 function~$f$, and to understand the deeper reasons that lead to the different types of
 series. This line of research was initiated by Bousquet-M\'{e}lou and Mishna~\cite{bousquet10}, and
 many other researchers have been contributing to it during the past years, see
 \cite{fayolle99,bousquet10,bostan10,bostan14b,kurkova15,bostan16b,courtiel17,bernardi17,dreyfus18} for
 some of the milestones. 

 As there are not many remaining open questions for the classical setting described above,
 people are now more and more turning to the study of variations and generalizations.
 One such generalization concerns the situation in higher dimensions.
 A first step was taken by Bostan and Kauers in~\cite{bostan09}, who used automated guessing
 to identify potentially D-finite step sets of size up to~5 in three dimensions. This work was extended by
 Bostan, Bousquet-Melou, Kauers, and Melczer~\cite{bostan16a} to step sets of size up to~6. They
 introduced the notion of a dimension of a lattice walk model, and the 
 so-called Hadamard decomposition of a step set, and they used these new concepts as well
 as the classical orbit sum method for proving the D-finiteness in certain cases.
 Bacher, Kauers and Yatchak~\cite{bacher16} have extended this work to step sets of arbitrary size,
 Du, Hou, and Wang provided non-D-finiteness results for many cases~\cite{du16}, 
 and most recently, Bogosel, Perrollaz, Raschel, and Trotignon~\cite{bogosel18} have systematically
 explored the asymptotic behaviour of counting sequences for walks in the octant and observed
 a striking relation between the nature of the generating function and the angles of certain
 triangles on the sphere.
 Despite all this progress, there are still many open questions related to walks
 in the octant. In particular, there is a list of 170 models whose nature remains unclear.
 For example, this list includes the 3D version of the classical 2D Kreweras model~\cite{kreweras65,bousquet05,bousquet10},
 the step set $\{(-1,0,0),(0,-1,0),(0,0,-1),(1,1,1)\}$. Although the 2D version has an
 algebraic generating function, the current asymptotic estimates suggest (without proof)
 that the 3D version is not D-finite.

 In this short note, we have nothing new to say about the 3D cases. Instead, our aim is to open
 the discussion for~4D. When the dimension of the lattice increases, the classification problem
 becomes more difficult in two ways. First, and most importantly, the total number of models
 explodes. For dimension~$D$, there are $2^{3^D-1}$ step sets, which evaluates to
 more than $10^{24}$ when $D=4$. There is no way to go through all of them in a reasonable
 time, even if we spend only a tiny amount of computation time per model. The second problem
 is that it won't be enough to spend only a tiny amount of computation time per model, because
 with increasing dimension it also becomes more costly to analyze a particular model. For example,
 computing the first $N$ terms of a counting sequence using the standard recurrence requires
 $\mathrm{O}(N^{D+1})$ time and $\mathrm{O}(N^D)$ memory. For $D=4$, this means that on a computer with 1\,Tb
 of main memory, we were only able to compute $N=700$ terms of a counting sequence. 

 
 \section{Search Procedure}\label{sec:search}

 In order to identify potentially interesting models, we have applied a similar search procedure
 as Bacher, Kauers, and Yatchak~\cite{bacher16} did in their search for interesting models in 3D. The
 procedure can be summarized as follows:
 \begin{itemize}
 \item \textbf{Only step sets $S\subseteq\{-1,0,1\}^4\setminus\{(0,0,0,0)\}$ with $|S|\leq7$ or $|S|\geq73$ were considered.}
   This restriction has no combinatorial motivation but was only made to reduce the computational
   cost to a manageable amount, similar as it was done in \cite{bostan09,bostan16a} for~3D.
 \item \textbf{Step sets containing unused steps were discarded.}
   Recall from~\cite{bostan16a} that an element $s$ of $S$ is called unused if it cannot appear in any walk of the model.
   For example, the step set $S=\{(1,0,-1,0),\penalty0(0,1,0,-1),\penalty0(1,1,0,0)\}$ leads to the same generating function
   as the step set $\{(1,1,0,0)\}$, because any use of $(1,0,-1,0)$ or $(0,1,0,-1)$ would lead the walk out
   of $\set N^4$, which is not allowed.
 \item \textbf{Only one step set from each symmetry class was considered.}
   Permuting the coordinates of the steps corresponds to a permutation of the variables of the generating
   function. For example, if $f(x_1,x_2,x_3,x_4,t)$ is the generating function for the model with step set
   $\{(1,0,1,1),\penalty0(-1,1,0,0),\penalty0(0,0,0,1)\}$, then $f(x_2,x_4,x_1,x_3,t)$ is the generating function for the
   model with step set $\{(0,1,1,1),\penalty0(1,0,-1,0),\penalty0(0,1,0,0)\}$. Since permutation of variables preserves
   algebraicity and D-finiteness, it suffices to consider one model per equivalence class. 
 \item \textbf{Step sets admitting a Hadamard decomposition were discarded.}
   Recall from~\cite{bostan16a} that a step set $S$ is said to admit a $(d_1,d_2)$-Hadamard decomposition for some
   positive $d_1,d_2$ with $d_1+d_2=D$ if it can be written as $S=(V\times\{0\})\cup(U\times W)$ with
   $V,U\subseteq\set Z^{d_1}$ and $W\subseteq\set Z^{d_2}$. If this is the case, the generating function
   for the lattice walk model for $S$ can be expressed in terms of the Hadamard product of the generating
   functions associated to the lower dimensional models corresponding to $U\cup V$ and $W$. 
 \item \textbf{Step sets with dimension less than 4 were discarded.}
   Recall from~\cite{bostan16a} that the dimension of a model is defined as the number of coordinates for which the
   nonnegativity restriction is not redundant. For example, for the step set $\{(1,1,1),\penalty0(1,-1,0),\penalty0(1,0,-1)\}$,
   the number of walks in $\set N^3$ is the same as the number of walks in $\set Z\times\set N^2$, because
   there is no way to get a negative first coordinate with the available steps. As the restriction on the other
   two coordinates is essential, the dimension is 2 in this case. Since lattice walk models in $\set N^4$
   whose dimension is less than 4 are equivalent to models in $\set N^3$ (possibly with multiple steps),
   it is fair to discard them. 
 \item \textbf{Step sets whose associated group has more than 800 elements were discarded.}
   Recall from~\cite{bousquet10,bostan16a} that to every model of maximal dimension we can associate a certain group.
   Given a step set $S\subseteq\{-1,0,1\}^D\setminus\{(0,\dots,0)\}$, the group is constructed as follows.
   For $i=1,\dots,D$, let $\Phi_i$ be the rational maps that sends $x_j$ to itself for $j\neq i$ and
   $x_i$ to $x_i^{-1}A_{i,-}/A_{i,+}$, where $A_{i,+}=\sum_{(s_1,\dots,s_D)\in S: s_i=1} x_1^{s_1}\cdots x_D^{s_D}/x_i$
   and $A_{i,-}=\sum_{(s_1,\dots,s_D)\in S: s_i=-1} x_1^{s_1}\cdots x_D^{s_D}/x_i^{-1}$.
   The group associated to $S$ is the group generated by $\Phi_1,\dots,\Phi_D$ under composition. A main
   result about the case $D=2$ is that this group is finite if and only if the generating function is D-finite~\cite{bousquet10,bostan16a,bostan16b}.
   While the experimental results for $D=3$ suggest that there may be non-D-finite cases with finite group,
   we are not aware of any (conjectured) D-finite case with an infinite group. For this reason, and also
   because a finite group gives the chance to apply the so-called orbit sum method for proving D-finiteness,
   we have decided to restrict the search to models with finite group. 
 \end{itemize}
 Out of the 7005847194 step sets with cardinality at most 7 or at least 73, there
 were 58 step sets which survived all these filters, the last filter being the by far strongest one.
 The surviving models are listed in the next section. They all have cardinality 5 or~7. 
 
 \section{Results}

 For models with a finite group, the orbit sum method is one approach to showing that the generating function
 is D-finite. It rests on the observation that, when certain technical conditions are satisfied, the generating
 function for a model can be expressed as
 \[
   f(x_1,\dots,x_D,t) = \frac1{x_1\cdots x_D}[x_1^>\dots x_D^>]\frac1{1-t P_S}\sum_{g\in G} g(x_1\cdots x_D),
 \]
 where $G$ is the group, $P_S:=\sum_{(s_1,\dots,s_D)\in S} x_1^{s_1}\cdots x_D^{s_D}$ is the step set polynomial
 (also called the inventory by some authors), and $[x_1^>\dots x_D^>]$ is the positive part extraction operator.
 Note that the expression to which the positive part extraction operator is applied is a rational function. By
 the closure of D-finiteness under taking positive parts, the formula above implies that the generating function
 is D-finite.

 For 50 of the 58 step sets identified by the procedure of Section~\ref{sec:search}, the orbit sum
 $\sum_{g\in G} g(x_1\cdots x_D)$ happens to be zero. In this case, the ``technical conditions'' alluded to above
 are not satisfied and we cannot directly conclude D-finiteness. In the other eight cases, we have checked
 with Yatchak's algorithm~\cite{yatchak17}
 that the technical conditions are satisfied, so the generating functions of these models are D-finite.

 For the 50 cases whose orbit sum is zero, we have tried to detect recurrence equations or differential equations
 via automated guessing, as systematically done in~\cite{bostan09} for 3D models. As remarked in the introduction, we
 were only able to compute 700 terms for each of these counting sequences, which only in one case
 (number~\ref{recurrence} in the listing below) was enough to find equations. For the generating function of walks
 with arbitrary endpoint,
 $f(1,\dots,1,t)$, we found a linear differential equation of order 12 with polynomial coefficients of degree up
 to~135. Its coefficient sequence appears to satisfy a linear recurrence of order 18 with polynomial coefficients
 of degree up to~113. 
 
 We suspect that further models are D-finite but only satisfy equations that are too large to be recovered from 700
 sequence terms, and we invite the lattice walk counting community to have a closer look at these models. In the
 tables below, we write $\bar 1$ instead of $-1$ for better readability. We also use a pictorial description of the
 step sets, extending similar descriptions used in the literature for lower dimensions. A step $(s_1,s_2,s_3,s_4)
 \in\{-1,0,1\}^4$ is represented by a bullet at position $(s_1,s_2,s_3,s_4)$, where $s_1$ is the column block
 ($-1=\mathrm{left}$, $0=\mathrm{middle}$, $1=\mathrm{right}$), $s_2$ is the row block ($1=\mathrm{top}$,
 $0=\mathrm{middle}$, $-1=\mathrm{bottom}$), and $s_3,s_4$ are the column and row, respectively, within the
 block specified by $s_1,s_2$. Models with nonzero orbit sum are highlighted. The orbit sums are stated in
 a separate table. 


\begin{tab}\em Models with a group isomorphic to $C_2\times C_2\times S_3$. \par\smallskip\noindent

   \global\let\nextcell\nextcella
   \begin{tabular}{@{}c|c|c|c@{}}%
     \hline
 \model{sAChAAAAAAEABCq}{---0,--0-,--++,-+00,+0--,+00+,+0+0}{C_2\times C_2\times S_3}{0}{1, 2, 6, 18, 84, 340}
 \model{sACMAAAAAAEAAYl}{---0,--0+,--+-,-+00,+0-+,+00-,+0+0}{C_2\times C_2\times S_3}{\frac{(w^2 - z) ( w z-1 ) (w - z^2) (w + w^2 z - w y^2 z + z^2) ( w^2 x^2 y - w^2 z + x^2 y z - w y^2 z - z^2 + w x^2 y z^2 - w )}{w^3 x y^2 z^3 (w^2 + z + w z^2)}}{1, 1, 3, 8, 33, 122}
 \model{sAEKAAAAAAEAAoQ}{----,--0+,--+0,-+00,+0-0,+00-,+0++}{C_2\times C_2\times S_3}{0}{1, 1, 4, 13, 58, 245}
 \model{sCAARgAAAAADEAM}{-0-0,-00+,-0+-,+--+,+-0-,++00,+-+0}{C_2\times C_2\times S_3}{\frac{(w^2 - z) (w z-1 ) (w - z^2) (w^2 + z - w y^2 z + w z^2) (w^2 x^2 - w y + x^2 z - w^2 y z + w x^2 y^2 z + w x^2 z^2 - y z^2)}{w^3 x y z^3 (w^2 + z + w y^2 z + w z^2)}}{ 1, 1, 3, 9, 35, 125}
 \model{sCAAUIAAAAAIUAD}{+---,-00-,-0-0,++00,+-+0,+-0+,-0++}{C_2\times C_2\times S_3}{0}{1, 1, 4, 14, 60, 238}
 \model{sCAAhQAAAAAFCAE}{-0--,+--0,+-0-,++00,-00+,-0+0,+-++}{C_2\times C_2\times S_3}{0}{1, 1, 4, 14, 63, 241}
 \model{sACMAAAAAAYgAEK}{--00,-+-0,-+0+,-++-,+0-+,+00-,+0+0}{C_2\times C_2\times S_3}{\frac{(w^2 - z) (1 - w z) (w - z^2) (w^2 x^2 y - w y^2 - w z + x^2 y z - w^2 y^2 z + w x^2 y z^2 - y^2 z^2) (w y^2 - w z + w^2 y^2 z +  y^2 z^2)}{w^2 x y^2 z^2 (w + w^2 z + z^2) (w^2 + z + w z^2)}}{1, 1, 4, 12, 62, 255}
 \model{sAChAAAAABCgAEE}{--00,-+-0,-+0-,-+++,+0--,+00+,+0+0}{C_2\times C_2\times S_3}{0}{1, 2, 8, 30, 166, 764}
 \model{sAEKAAAAAAoQAEJ}{-+--,--00,-+0+,-++0,+0-0,+00-,+0++}{C_2\times C_2\times S_3}{0}{1, 1, 6, 21, 126, 581}
 \model{sUIACAAAAAAIUAF}{-0-0,-00-,-0++,++--,+-00,++0+,+++0}{C_2\times C_2\times S_3}{0}{\scalebox{.86}[1.0]{1, 2, 10, 46, 260, 1402}}
 \model{shQACAAAAAAFCAH}{-0--,-00+,-0+0,++-0,++0-,+-00,++++}{C_2\times C_2\times S_3}{0}{1, 1, 7, 33, 197, 1065}
 \model{sRgACAAAAAADEAG}{-0-0,-00+,-0+-,++-+,+-00,++0-,+++0}{C_2\times C_2\times S_3}{\frac{(w^2 - z) (1 - w z) (w - z^2) (w^2 y^2 - w z + y^2 z + w y^2 z^2) (w^2 x^2 y^2 + w x^2 z - w^2 y z + x^2 y^2 z - y z^2 + w x^2 y^2 z^2-w y )}{w^2 x y z^2 (w^2 + z + w z^2) (w^2 y^2 + w z + y^2 z + w y^2 z^2)}}{1, 1, 5, 20, 102, 496}
 \end{tabular}
 \end{tab}
 
\begin{tab}\em Models with a group isomorphic to $S_3\times S_3$. \par\smallskip\noindent

   \global\let\nextcell\nextcella
   \begin{tabular}{@{}c|c|c|c@{}}%
     \hline
 
 \model{recurrence}{-0-0,0--+,-00+,0-0-,++00,-0+-,0-+0}{S_3\times S_3}{0}{1, 1, 3, 9, 27, 117}
 \model{sCAAAAAAIUAFCAD}{-0--,0--0,0-0-,++00,-00+,-0+0,0-++}{S_3\times S_3}{0}{1, 1, 4, 14, 45, 223}
 \model{sAAQAChAAAAABCr}{---0,--0-,--++,0+--,0+0+,0++0,+000}{S_3\times S_3}{0}{1, 3, 9, 37, 169, 759}
 \model{sAAQACMAAAAAAYl}{---0,--0+,--+-,0+-+,0+0-,0++0,+000}{S_3\times S_3}{0}{1, 2, 5, 18, 72, 295}
 \model{sAAQAEKAAAAAAoW}{----,--0+,--+0,0+-0,0+0-,0+++,+000}{S_3\times S_3}{0}{1, 2, 6, 26, 118, 548}
 \model{sAAARgQAAAADEAF-flipped}{0--0,0-0+,0-+-,+000,-+-+,-+0-,-++0}{S_3\times S_3}{\frac{(w^2 - z) (1 - w z) (w - z^2) (w + w^2 z - w x y z + z^2) (w^2 x + x z - w y^2 z + w x z^2) (w^2 x^2 - w y + x^2 z - w^2 y z + w x^2 z^2 - y z^2)}{w^4 x^2 y^2 z^4 (w^2 + z + w z^2)}}{ 1, 1, 2, 5, 14, 47}
 \model{sAAQAAAAIUoQAAD}{-+--,0--0,0-0-,+000,-+0+,-++0,0-++}{S_3\times S_3}{0}{1, 1, 3, 9, 27, 103}
 \model{sAAQAAAAFDCgAAH}{0---,-+0-,-+-0,+000,0-+0,0-0+,-+++}{S_3\times S_3}{0}{1, 1, 2, 6, 19, 73}
 \model{sACMABiAAAAAAEC}{--00,0+-0,0+0+,0++-,+0-+,+00-,+0+0}{S_3\times S_3}{0}{1, 2, 8, 36, 184, 978}
 \model{sAEKAChAAAAAAEH}{--00,0+--,0+0+,0++0,+0-0,+00-,+0++}{S_3\times S_3}{0}{\scalebox{0.86}[1.0]{1, 3, 14, 74, 425, 2515}}
 \model{sAChAAAAAhCgAAM}{-+-0,-+0-,-+++,0-00,+0--,+00+,+0+0}{S_3\times S_3}{0}{1, 2, 8, 34, 176, 908}
 \model{sAEKAAAAAgoQAAE}{-+--,+0-0,+00-,0-00,-+0+,-++0,+0++}{S_3\times S_3}{0}{1, 1, 6, 24, 133, 695}
 \model{sACMAAAAAgYgAAF}{-+-0,-+0+,-++-,0-00,+0-+,+00-,+0+0}{S_3\times S_3}{\frac{(w^2\! - z) (1\! - w z) (w\! - z^2) (w^2 x^2\! - w y + x^2 z - w^2 y z + w x^2 z^2\! - y z^2) (w^2 x y - w z + x y z + w x y z^2) (w y^2 - w x z + w^2 y^2 z + y^2 z^2)}{w^2 x^2 y^2 z^2 (w + w^2 z +  z^2) (w^2 + z + w z^2)^2}}{ 1, 1, 4, 14, 66, 309}
 \model{sUIAAAAAAgAIUAL}{-0-0,-00-,-0++,0-00,++--,++0+,+++0}{S_3\times S_3}{0}{\scalebox{0.86}[1.0]{1, 2, 10, 46, 244, 1358}}
 \model{shQAAAAAAgAFCAA}{-0--,-00+,-0+0,0-00,++-0,++0-,++++}{S_3\times S_3}{0}{1, 1, 7, 33, 181, 1025}
 \model{sRgAAAAAAgADEAM}{-0-0,-00+,-0+-,0-00,++-+,++0-,+++0}{S_3\times S_3}{0}{1, 1, 5, 20, 94, 478}
 \model{sAAACEKAAAAFCAG}{-0--,0+-0,0+0-,+-00,-00+,-0+0,0+++}{S_3\times S_3}{0}{1, 1, 4, 16, 65, 299}
 \model{sAAACChAAAAIUAC}{-0-0,-00-,-0++,0+--,0+0+,0++0,+-00}{S_3\times S_3}{0}{1, 2, 6, 22, 94, 414}
 \model{sAAACCMAAAADEAA}{-0-0,-00+,-0+-,0+-+,0+0-,0++0,+-00}{S_3\times S_3}{\frac{(w^2 - z) ( w z-1 ) (w - z^2) (w y - w x^2 z + w^2 y z + y z^2) ( w^2 x y - w^2 z + x y z - z^2 + w x y z^2-w ) (w^2 y^2 - w x z + y^2 z + w y^2 z^2)}{w^3 x^2 y^2 z^3 (w^2 + z + w z^2)^2}}{1, 1, 3, 10, 37, 151}

 \\\hline
 
 \end{tabular}
 \end{tab}

 \begin{tab}\em Models with a group isomorphic to $S_5$. \par\smallskip\noindent

   \global\let\nextcell\nextcella
   \begin{tabular}{@{}c|c|c|c@{}}%
     \hline

 \model{sAAQAAASAgBAAAA}{-+-+,0-00,000-,00+0,+000}{S_5}{0}{ 1, 2, 4, 10, 30, 98}
 \model{sgAAAAACggAAgAI}{-000,0-00,000-,00-0,++++}{S_5}{0}{ 1, 1, 5, 21, 81, 325}
 \model{sAAQAAQUAAAAAAT}{----,00+0,000+,0+00,+000}{S_5}{0}{\scalebox{0.86}[1.0]{1, 4, 16, 64, 256, 1048}}
 \model{sAAQAAAUAgAQAAL}{-+--,0-00,000+,00+0,+000}{S_5}{0}{ 1, 3, 9, 27, 87, 303}
 \model{sAAQAAACghAAAAL}{0-00,00-0,000-,+000,-+++}{S_5}{0}{ 1, 1, 2, 6, 21, 73}


 \model{sAACAAASBAEAAAI}{000-,+0-0,-+00,00+0,0-0+}{S_5}{\frac{ (w^2 - y) ( w y-x ) (w x - y^2) (w x - z) (  w z-1) ( w z-x y) ( x z-y) (w - y z) (x^2 - y z) (x - z^2)}{w^4 x^4 y^4 z^4}}{1, 1, 2, 4, 10, 26}
 \model{sCAAAAAChAAEAAE}{000-,00-0,++00,-0+0,0-0+}{S_5}{0}{ 1, 1, 3, 9, 29, 99}
 \model{sAACAAIUAAAAAEH}{--00,00+0,000+,0+0-,+0-0}{S_5}{0}{ 1, 2, 6, 18, 60, 206}
 \model{sACAAAAEgQEAAAG}{0-0-,00-0,+0+0,-+00,000+}{S_5}{0}{ 1, 2, 6, 20, 71, 269}
 \model{sCAAAAAUAQAAEAG}{-0-0,0-0-,000+,00+0,++00}{S_5}{0}{ 1, 3, 9, 31, 117, 467}
 \model{sACAAAgCgAAAAEH}{--00,000-,00-0,+0+0,0+0+}{S_5}{0}{ 1, 2, 6, 20, 80, 318}
 \model{sACAAAAChAEAAAI}{00-0,000-,-+00,0-0+,+0+0}{S_5}{0}{ 1, 1, 3, 8, 24, 78}
 \model{sAACAAAUAQEAAAG}{0-0-,+0-0,00+0,-+00,000+}{S_5}{0}{ 1, 2, 5, 14, 42, 136}
 \model{sCAAAAASBAAAEAI}{-0-0,000-,0-0+,00+0,++00}{S_5}{0}{  1, 2, 5, 16, 57, 209}
 \model{sACAAAIEgAAAAEH}{--00,00-0,000+,0+0-,+0+0}{S_5}{0}{ 1, 2, 6, 18, 63, 229}
		
		
 \model{sAAQAAASBAAgAAA}{-+-0,000-,0-0+,00+0,+000}{S_5}{0}{1, 2, 4, 10, 28, 82}
 \model{sAACAAASAgIAAAH}{0-00,000-,+0-0,00+0,-+0+}{S_5}{0}{ 1, 1, 2, 4, 11, 31}
 \model{sQAAAAAChAAAgAH}{-000,00-0,000-,0-0+,+++0}{S_5}{0}{ 1, 1, 4, 14, 49, 183}
 \model{sAACAAQUAAAAACG}{--0-,00+0,+0-0,000+,0+00}{S_5}{0}{ 1, 3, 10, 35, 126, 474}
 \model{sACAAAAEggCAAAA}{-+0-,00-0,0-00,000+,+0+0}{S_5}{0}{ 1, 2, 5, 15, 52, 185}
 \model{sAAQAAAEgQgAAAF}{0-0-,00-0,+000,000+,-++0}{S_5}{0}{ 1, 2, 5, 14, 45, 159}
 \model{sAAQAAAUAQAgAAL}{0-0-,00+0,+000,-+-0,000+}{S_5}{0}{ 1, 3, 9, 29, 99, 355}
 \model{sAACAAAUAgCAAAA}{-+0-,0-00,00+0,+0-0,000+}{S_5}{0}{ 1, 2, 5, 13, 38, 119}
 \model{sACAAAACggIAAAH}{0-00,000-,00-0,+0+0,-+0+}{S_5}{0}{ 1, 1, 3, 9, 31, 109}
 \model{sAAQAAAChAgAAAH}{00-0,000-,+000,0-0+,-++0}{S_5}{0}{ 1, 1, 2, 5, 15, 47}
 \model{sEAAAAASAgAAEAF}{-0-0,0-00,000-,00+0,++0+}{S_5}{0}{ 1, 2, 6, 22, 88, 358}
 \model{sACAAAQEgAAAACG}{--0-,00-0,000+,0+00,+0+0}{S_5}{0}{ 1, 3, 10, 35, 132, 534}

 \\\hline
 
 \end{tabular}
 \end{tab}

\begin{tab}\em Nonzero orbit sums. \par\smallskip\noindent
  \begin{tabular}{c|c}
  idx & orbit sum \\\hline
 \Model{sCAARgAAAAADEAM}{-0-0,-00+,-0+-,+--+,+-0-,++00,+-+0}{C_2\times C_2\times S_3}{\frac{(w^2 - z) (w z-1 ) (w - z^2) (w^2 + z - w y^2 z + w z^2) (w^2 x^2 - w y + x^2 z - w^2 y z + w x^2 y^2 z + w x^2 z^2 - y z^2)}{w^3 x y z^3 (w^2 + z + w y^2 z + w z^2)}}{ 1, 1, 3, 9, 35}
 \Model{sACMAAAAAAEAAYl}{---0,--0+,--+-,-+00,+0-+,+00-,+0+0}{C_2\times C_2\times S_3}{\frac{(w^2 - z) ( w z-1 ) (w - z^2) (w + w^2 z - w y^2 z + z^2) ( w^2 x^2 y - w^2 z + x^2 y z - w y^2 z - z^2 + w x^2 y z^2 - w )}{w^3 x y^2 z^3 (w^2 + z + w z^2)}}{1, 1, 3, 8, 33}
 \Model{sACMAAAAAAYgAEK}{--00,-+-0,-+0+,-++-,+0-+,+00-,+0+0}{C_2\times C_2\times S_3}{\frac{(w^2 - z) (1 - w z) (w - z^2) (w^2 x^2 y - w y^2 - w z + x^2 y z - w^2 y^2 z + w x^2 y z^2 - y^2 z^2) (w y^2 - w z + w^2 y^2 z +  y^2 z^2)}{w^2 x y^2 z^2 (w + w^2 z + z^2) (w^2 + z + w z^2)}}{1, 1, 4, 12, 62}
 \Model{sRgACAAAAAADEAG}{-0-0,-00+,-0+-,++-+,+-00,++0-,+++0}{C_2\times C_2\times S_3}{\frac{(w^2 - z) (1 - w z) (w - z^2) (w^2 y^2 - w z + y^2 z + w y^2 z^2) (w^2 x^2 y^2 + w x^2 z - w^2 y z + x^2 y^2 z - y z^2 + w x^2 y^2 z^2-w y )}{w^2 x y z^2 (w^2 + z + w z^2) (w^2 y^2 + w z + y^2 z + w y^2 z^2)}}{1, 1, 5, 20, 102}
 \Model{sAAARgQAAAADEAF-flipped}{0--0,0-0+,0-+-,+000,-+-+,-+0-,-++0}{S_3\times S_3}{\frac{(w^2 - z) (1 - w z) (w - z^2) (w + w^2 z - w x y z + z^2) (w^2 x + x z - w y^2 z + w x z^2) (w^2 x^2 - w y + x^2 z - w^2 y z + w x^2 z^2 - y z^2)}{w^4 x^2 y^2 z^4 (w^2 + z + w z^2)}}{ 1, 1, 2, 5, 14, 47}
 \Model{sACMAAAAAgYgAAF}{-+-0,-+0+,-++-,0-00,+0-+,+00-,+0+0}{S_3\times S_3}{\frac{(w^2\! - z) (1\! - w z) (w\! - z^2) (w^2 x^2\! - w y + x^2 z - w^2 y z + w x^2 z^2\! - y z^2) (w^2 x y - w z + x y z + w x y z^2) (w y^2 - w x z + w^2 y^2 z + y^2 z^2)}{w^2 x^2 y^2 z^2 (w + w^2 z +  z^2) (w^2 + z + w z^2)^2}}{ 1, 1, 4, 14, 66}
 \Model{sAAACCMAAAADEAA}{-0-0,-00+,-0+-,0+-+,0+0-,0++0,+-00}{S_3\times S_3}{\frac{(w^2 - z) ( w z-1 ) (w - z^2) (w y - w x^2 z + w^2 y z + y z^2) ( w^2 x y - w^2 z + x y z - z^2 + w x y z^2-w ) (w^2 y^2 - w x z + y^2 z + w y^2 z^2)}{w^3 x^2 y^2 z^3 (w^2 + z + w z^2)^2}}{1, 1, 3, 10, 37}
 \Model{sAACAAASBAEAAAI}{000-,+0-0,-+00,00+0,0-0+}{S_5}{\frac{ (w^2 - y) ( w y-x ) (w x - y^2) (w x - z) (  w z-1) ( w z-x y) ( x z-y) (w - y z) (x^2 - y z) (x - z^2)}{w^4 x^4 y^4 z^4}}{1, 1, 2, 4, 10, 26}
  \end{tabular} 
 \end{tab}


 \bibliographystyle{plain}
 \bibliography{bib}

\begin{thebibliography}{10}

\bibitem{bacher16}
Axel Bacher, Manuel Kauers, and Rika Yatchak.
\newblock Continued classification of 3d lattice walks in the positive octant.
\newblock In {\em Proceedings of FPSAC'16}, pages 95--105, 2016.

\bibitem{bernardi17}
Olivier Bernardi, Mireille Bousquet-M{\'e}lou, and Kilian Raschel.
\newblock Counting quadrant walks via {T}utte's invariant method.
\newblock Technical Report 1708.08215, ArXiv, 2017.

\bibitem{bogosel18}
Beniamin Bogosel, Vincent Perrollaz, Kilian Raschel, and Am{\'e}lie Trotignon.
\newblock 3d positive lattice walks and spherical triangles.
\newblock Technical Report 1807.08610, ArXiv, 2018.

\bibitem{bostan16a}
Alin Bostan, Mireille Bousquet-M{\'e}lou, Manuel Kauers, and Stephen Melczer.
\newblock On 3-dimensional lattice walks confined to the positive octant.
\newblock {\em Annals of Combinatorics}, 20(4):661--704, 2016.

\bibitem{bostan16b}
Alin Bostan, {Fr\'ed\'eric} Chyzak, Mark van Hoeij, Manuel Kauers, and Lucien
  Pech.
\newblock Hypergeometric expressions for generating functions of walks with
  small steps in the quarter plane.
\newblock {\em European Journal of Combinatorics}, 61:242--275, 2017.

\bibitem{bostan09}
Alin Bostan and Manuel Kauers.
\newblock Automatic classification of restricted lattice walks.
\newblock In {\em Proceedings of FPSAC'09}, pages 201--215, 2009.

\bibitem{bostan14b}
Alin Bostan, Kilian Raschel, and Bruno Salvy.
\newblock Non-d-finite excursions in the quarter plane.
\newblock {\em Journal of Combinatorial Theory, Series A}, 121:45--63, 2014.

\bibitem{bostan10}
Alin Bostan and Manuel~Kauers with an appendix by Mark~van Hoeij.
\newblock The complete generating function for {G}essel walks is algebraic.
\newblock {\em Proceedings of the AMS}, 138(9):3063--3078, 2010.

\bibitem{bousquet05}
Mireille Bousquet-Melou.
\newblock Walks in the quarter plane: {K}reweras' algebraic model.
\newblock {\em The Annals of Applied Probability}, 15(2):1451--1491, 2005.

\bibitem{bousquet10}
Mireille Bousquet-M{\'e}lou and Marni Mishna.
\newblock Walks with small steps in the quarter plane.
\newblock {\em Contemporary Mathematics}, 520:1--40, 2010.

\bibitem{courtiel17}
Julien Courtiel, Stephen Melczer, Marni Mishna, and Kilian Raschel.
\newblock Weighted lattice walks and universality classes.
\newblock {\em Journal of Combinatorial Theory, Series A}, 152:255--302, 2017.

\bibitem{dreyfus18}
Thomas Dreyfus, Charlotte Hardoin, Julien Roques, and Michael~F. Singer.
\newblock On the nature of the generating series of walks in the quarter plane.
\newblock {\em Inventiones mathematicae}, 213(1):205--236, 2018.

\bibitem{du16}
D.K. Du, Q.-H. Hou, and R.-H. Wang.
\newblock Infinite orders and non-d-finite property of 3-dimensional lattice
  walks.
\newblock {\em Electronic Journal of Combinatorics}, 23, 2016.

\bibitem{fayolle99}
Guy Fayolle, Roudolf Iasnogorodski, and Vadim Malyshev.
\newblock {\em Random Walks in the quarter-plane}.
\newblock Springer, 1999.

\bibitem{kreweras65}
Germain Kreweras.
\newblock Sur une classe de probl{\`e}mes li{\'e}s au treillis des partitions
  d'entiers.
\newblock {\em Cahiers du B.U.R.O.}, 6:5--105, 1965.

\bibitem{kurkova15}
Irina Kurkova and Kilian Raschel.
\newblock New steps in walks with small steps in the quarter plane.
\newblock {\em Annals of Combinatorics}, 19:461--511, 2015.

\bibitem{yatchak17}
Rika Yatchak.
\newblock {\em Restricted Lattice Paths and Their Generating Functions}.
\newblock PhD thesis, Johannes Kepler University Linz, 2017.

\end{thebibliography}

 \end{document}